\documentclass[11pt, twoside]{article}
\usepackage{amssymb}
\usepackage{amsmath, amssymb}\usepackage{cite}\usepackage{mathrsfs}\usepackage[amsmath, thmmarks]{ntheorem}
\usepackage{listings}
\usepackage[titletoc]{appendix}\allowdisplaybreaks
\usepackage{multirow}

\makeatletter
\renewcommand\theequation{\thesection.\@arabic\c@equation}
\makeatother

\newtheorem{thm}{Theorem}[section]%
\newtheorem{lem}[thm]{Lemma}%
\input cyracc.def

\def\f{\noindent}

\def\demo{\f{\bf Proof}\hskip10pt}

\def\qed{\hfill $\Box$}

\begin{document}

\title{{\bf Cuntz algebra automorphisms: transpositions}}
\footnotetext{E-mail addresses: Junyao$_{-}$Pan@126.com}

\author{{\bf Junyao Pan}\\
{\footnotesize School of Cyber Science and Engineering, Wuxi University, Wuxi, Jiangsu, 214105, P. R. China}
}

\date{}
\maketitle

\noindent{\small {\bf Abstract:} Permutative automorphisms of the Cuntz algebras $\mathcal{O}_n$ are in bijection with the stable permutations of $[n]^t$. They are also the elements of the reduced Weyl group of $Aut(\mathcal{O}_n)$. In this paper, we characterize the stability of transpositions in $S([n]^3)$, and thus providing a new family (with $6$ degrees of freedom) of automorphisms of the Cuntz algebras $\mathcal{O}_n$ for any $n>1$.

\vskip0.2cm
\noindent{\small {\bf Keywords:} Cuntz algebra; Permutative automorphism; Stable permutation.

\vskip0.2cm
\noindent{\small {\bf Mathematics Subject Classification (2020):} 05E16, 05A05, 05A15}

\section {Introduction}

Throughout this paper, $N$ denotes the set of all positive integers; $[n]:=\{1,2,...,n\}$; $[n]^{t}$ stands for the cartesian product of $t$ copies of $[n]$ for $t\in N$; $S([n]^{t})$ expresses the symmetric group on $[n]^{t}$.

A $C^*$-algebra is a norm closed self-adjoint sub-algebra of the bounded operators on a Hilbert space. Due to a theorem of Gelfand, Naimark and Segal, one can alternatively describe $C^*$-algebras axiomatically as complex Banach algebras with an involution that satisfies $\|a^*a\|=\|a\|^2$. Actually, $C^*$-algebras were first introduced for providing a suitable environment for a rigorous approach to quantum theories \cite{HAA}, and more recently have been related to such diverse fields as operator theory, group representations, topology, quantum mechanics, non-commutative geometry and dynamical systems \cite{CON}. On the other hand, the symmetries of $C^*$-algebras are provided by automorphisms in the most classical sense. Despite the fact that many general results and constructions about automorphisms of $C^*$-algebras have been obtained in \cite{FA,PE} and that they are used extensively, not very much is known about the automorphism group of most $C^*$-algebras. In addition, many other interesting (simple) examples of $C^*$-algebras emerged in the 1970s and 1980s, such as Cuntz \cite{C1} invented his Cuntz algebras $\mathcal{O}_n$. These are, for $2 \leq n < \infty$, generated by $n$ isometries $s_1, s_2,..., s_n$ satisfying the relation $I = s_1s^*_1+s_2s^*_2+\cdot\cdot\cdot+s_ns^*_n$. For $n = \infty$, this relation is replaced with the relation that the support projections $s_js^*_j$ are mutually orthogonal. Actually, they have been attracting much attention since their appearance in the seminal paper, such as \cite{BJ,C2,C3,CK,DR1,DR2}. Moreover, the study of automorphisms of the Cuntz algebras soon appeared quite intriguing and revealed many interesting facets. Notably, following the deep insight by Cuntz \cite{C4}, in a series of papers \cite{CHS1,CHS2,CS,CKS} a general theory of restricted Weyl groups for $Aut(\mathcal{O}_n)$ has been studied both from a theoretical viewpoint as well as from the perspective of constructing explicit examples. In order to explain the connection with the main topic of this note, we next introduce some more terminology.

 Let $\mathcal{O}_n$ be the Cuntz algebra with $n \geq2$, $\mathcal{F}_n \subset\mathcal{O}_n$ the so-called core-UHF $C^*$-subalgebra generated by a nested family of subalgebras $\mathcal{F}^k_n$ that is isomorphic to the algebra of complex matrices $M_{n^k}$ and $\mathcal{D}_n$ be the $C^*$-subalgebra generated by the family of subalgebras $\mathcal{D}^k_n$ of $\mathcal{F}^k_n$ that is isomorphic to the algebra of diagonal matrices in $M_{n^k}$. Then, following the insight by Cuntz, the reduced Weyl group of $Aut(\mathcal{O}_n)$ is defined as $$Aut(\mathcal{O}_n,\mathcal{F}_n)\cap Aut(\mathcal{O}_n,\mathcal{D}_n)/Aut_{\mathcal{D}_n}(\mathcal{O}_n),$$ where $Aut(\mathcal{O}_n, X)$ is the subgroup of $Aut(\mathcal{O}_n)$ consisting of the automorphisms which leave $X$ invariant, while $Aut_X(\mathcal{O}_n)$ is the one of those which fix $X$ pointwise. It is well-known that unital $*$-endomorphisms of $\mathcal{O}_n$ are in bijective correspondence with unitaries in $\mathcal{O}_n$, call this bijection $u\mapsto\lambda_u$. With some more work the reduced Weyl group can be further identified with the set of automorphisms $\lambda_u$ of $\mathcal{O}_n$ induced by the so-called permutative unitaries $u \in \cup_{k\geq1}\mathcal{F}^k_n$. Moreover, for general permutative unitaries $u \in \mathcal{F}^k_n$, for any $k$, it was shown in \cite{CHS2}, that $\lambda_u$ is an automorphism precisely when the sequence of unitaries $(\varphi^r(u^*)\cdot\cdot\cdot\varphi(u^*)u^*\varphi(u)\cdot\cdot\cdot\varphi^r(u))_{r\geq0}$ in $\mathcal{F}_n$ eventually stabilizes, where the endomorphism $\varphi$ of $\mathcal{F}_n$ corresponds, in the isomorphism between $M_{n^k}$ and $M_n\otimes M_n\otimes\cdot\cdot\cdot\otimes M_n$ with $k$ factors, to the tensor shift map $x \mapsto 1_{M_n}\otimes x$. By coherently identifying permutative unitaries in $\mathcal{F}^k_n$ with permutation matrices in $M_{n^k}$ and thus with permutations of the set $\{1, ..., n^k\}$, and finally with permutations of the set $[n]^k$ by lexicographic ordering, one is lead to the class of \emph{stable permutations} that were defined by Brenti and Conti \cite{BC}. Notably, Brenti and Conti \cite{BC} pointed out that $\lambda_u$ is an automorphism of $\mathcal{O}_n$ if and only if $u$ is stable, in particular, the reduced Weyl group of $\mathcal{O}_n$ is $\{\lambda_u : u \in S([n]^r), r \in N, u~{\rm{ stable}}\}$. This opened up the investigations of the reduced Weyl groups of the Cuntz algebras $\mathcal{O}_n$ from a combinatorial point of view. Following the deep insight, Brenti and Conti and Nenashev continued to deep these investigations and focus on the explicit construction of restricted Weyl group elements using combinatorial techniques in \cite{BCN1,BCN2}. Next we recall some notions and notations about stable permutations, for more detailed information see \cite{BC}.

Given two permutations $u\in S([n]^t)$ and $v\in S([n]^r)$ with $n\in N$ and $n\geq2$, the \emph{tensor product} of $u$ and $v$ is the permutation $u\otimes v\in S([n]^{t+r})$ defined by $$(u\otimes v)(\alpha,\beta):=(u(\alpha),v(\beta))$$ for all $\alpha\in[n]^t$ and $\beta\in[n]^r$. For a permutation $u \in S([n]^t)$, define a sequence of permutations $\Psi_k(u) \in S([n]^{t+k})$, $k \geq 0$ by setting $\Psi_0(u) := u^{-1}$ and, for $k \in N$, $$\Psi_k(u)=\prod_{i = 0}^{k}\underbrace{1\otimes\cdot\cdot\cdot\otimes1}_{k-i}\otimes u^{-1}\otimes\underbrace{1\otimes\cdot\cdot\cdot\otimes1}_{i}\prod_{i = 1}^{k}\underbrace{1\otimes\cdot\cdot\cdot\otimes1}_{i}\otimes u\otimes\underbrace{1\otimes\cdot\cdot\cdot\otimes1}_{k-i},$$ where $1$ denotes the identity of $S_n := S([n])$. Then $u$ as above is said to be \emph{stable} if there exists some integer $k \geq 1$ such that $$\Psi_{k+l}(u)=\Psi_{k-1}(u)\otimes\underbrace{1\otimes\cdot\cdot\cdot\otimes1}_{l+1}~{\rm{for ~all}} ~l\geq0,$$
and the least such value of $k$ is called the \emph{rank} of $u$, denoted by $rank(u)$. One easily checks that if $t=1$ then all permutations in $S_n$ are stable of rank $1$. However, it seems to be extremely difficult to determine whether $u$ is stable for $t>1$, and even the stability of $r$-cycles in $S([n]^2)$ is still open for $r\geq6$. Actually, the current research mainly focuses on the the stability of special permutations in $S([n]^2)$, such as Brenti and Conti and Nenashev \cite{BCN1,BCN2} investigated the stable $r$-cycles in $S([n]^2)$ for $r\leq5$, and the author \cite{PA} gave a description of the stability of some involutions in $S([n]^2)$. Thereby, Brenti and Conti and Nenashev \cite{BCN2} proposed an open problem that examine the stability of cycles in $S([n]^t)$ for $t>2$.

In this paper, we investigate the stability of $2$-cycles in $S([n]^3)$ and completely characterize the stability of transpositions in $S([n]^3)$, and thus providing a new family (with $6$ degrees of freedom) of automorphisms of the Cuntz algebras $\mathcal{O}_n$ for any $n>1$. The result is as follows.

\begin{thm}\label{pan1-2}\normalfont
Let $u=((a_1,a_2,a_3),(b_1,b_2,b_3))$ be a transposition in $S([n]^3)$. Then the following conditions are equivalent:\\\\
${\bf{(1)}}$ $u$ is stable;\\
${\bf{(2)}}$ $u$ is stable of rank $1$;\\
${\bf{(3)}}$ either $\{a_1, b_1\} \cap \{a_3, b_3\}= \emptyset$ and $\{a_1, a_2\}\neq\{b_2, b_3\}$ and $\{b_1, b_2\}\neq\{a_2, a_3\}$ or\\
$a_1=a_3=b_1=b_3$ and $a_1\neq a_2\neq b_2\neq b_1$.
\end{thm}

\section {Proof of Theorem 1.1}

First of all we present some simple facts that will be used repeatedly. In particular, the first one provides a more specific characterization of the stable permutations of rank $1$ in $S([n]^3)$, which can be seen as an extension of \cite[Proposition 4.5]{BC} on $S([n]^3)$.

\begin{lem}\label{pan2-1}\normalfont
Given a permutation $u$ in $S([n]^3)$. Then $u$ is stable of rank $1$ if and only if $u$ satisfies the equations
\begin{equation}
\begin{cases}
\Psi_{1}(u)=\Psi_{0}(u)\otimes1 \\
\Psi_{2}(u)=\Psi_{1}(u)\otimes1
\end{cases}
\end{equation}
\end{lem}
\demo Apply the definitions of the stable permutation and its rank, this lemma can be deduced immediately.  \qed

\begin{lem}\label{pan2-11}\normalfont
Let $u=((a_1,a_2,a_3),(b_1,b_2,b_3))$ be a transposition in $S([n]^3)$. Then
\begin{equation*}
\begin{cases}
u^{-1}\otimes1=u\otimes1=\prod_{x=1}^{n}((a_1,a_2,a_3,x),(b_1,b_2,b_3,x)) \\\\
u^{-1}\otimes1\otimes1=u\otimes1\otimes1=\prod_{x=1}^{n}\prod_{y=1}^{n}((a_1,a_2,a_3,x,y),(b_1,b_2,b_3,x,y))
\end{cases}
\end{equation*}
\end{lem}
\demo By definition of the tensor product, we deduce this lemma immediately.  \qed

\begin{lem}\label{pan2-111}\normalfont([2, Proposition 4.4])
Let $u \in S([n]^{t})$, and $k \in N$. Then\\\\
$i)$ if $\Psi_{k}(u)\in S([n]^{k+1})\otimes\underbrace{1\otimes\cdot\cdot\cdot\otimes1}_{t-1}$ then $u$ is stable of rank $\leq k + 1$;\\
$ii)$ if $u$ is stable of rank $\leq k + 1$ then $\Psi_{k+1}(u)\in S([n]^{k+t})\otimes1$.\\\\
In particular, $u$ is stable if and only if there exists a positive integer $h$ such that $$\Psi_{h}(u)\in S([n]^{h+1})\otimes\underbrace{1\otimes\cdot\cdot\cdot\otimes1}_{t-1}.$$
\end{lem}

Now we start to prove Theorem\ \ref{pan1-2} by ${\bf{(3)}}$ implies ${\bf{(2)}}$.

\begin{lem}\label{pan2-2}\normalfont
Let $u=((a_1,a_2,a_3),(b_1,b_2,b_3))$ be a $2$-cycle in $S([n]^3)$ with $\{a_1, b_1\} \cap \{a_3, b_3\}= \emptyset$ and $\{a_1, a_2\}\neq\{b_2, b_3\}$ and $\{b_1, b_2\}\neq\{a_2, a_3\}$. Then $u$ is stable of rank $1$.
\end{lem}
\demo By Lemma\ \ref{pan2-1}, it suffices to confirm that $u$ satisfies the equations (2.1). According to Lemma\ \ref{pan2-11} and definition of $\Psi_{k}(u)$, we deduce that $$\Psi_{1}(u)=(1\otimes u^{-1})(u^{-1}\otimes 1)(1\otimes u)=(1\otimes u^{-1})\Big(\prod_{x=1}^{n}((a_1,a_2,a_3,x),(b_1,b_2,b_3,x))\Big)(1\otimes u).$$
Since $1\otimes u^{-1}=(1\otimes u)^{-1}$, it follows that
$$(1\otimes u^{-1})(u^{-1}\otimes 1)(1\otimes u)=\prod_{x=1}^{n}((a_1,u(a_2,a_3,x)),(b_1,u(b_2,b_3,x))).$$
Note that $u(a_2,a_3,x)\neq(a_2,a_3,x)$ if and only if $(a_2,a_3,x)=(a_1,a_2,a_3)$ or $(a_2,a_3,x)=(b_1,b_2,b_3)$. If $(a_2,a_3,x)=(a_1,a_2,a_3)$ then $a_1=a_2=a_3$. In this case, $a_1\in\{a_1, b_1\} \cap \{a_3, b_3\}$, a contradiction. If $(a_2,a_3,x)=(b_1,b_2,b_3)$ then $a_2=b_1$ and $a_3=b_2$, and thus $\{b_1, b_2\}=\{a_2, a_3\}$, a contradiction. Thereby, $u(a_2,a_3,x)=(a_2,a_3,x)$ for all $x\in[n]$. Similarly, $u(b_2,b_3,x)=(b_2,b_3,x)$ for all $x\in[n]$. Hence, we have
$$\Psi_{1}(u)=(1\otimes u^{-1})(u^{-1}\otimes 1)(1\otimes u)=\prod_{x=1}^{n}((a_1,a_2,a_3,x),(b_1,b_2,b_3,x))=u^{-1}\otimes 1=\Psi_{0}(u)\otimes1.$$

Consider $\Psi_{2}(u)=(1\otimes1\otimes u^{-1})(1\otimes u^{-1}\otimes1)(u^{-1}\otimes1\otimes1)(1\otimes u\otimes1)(1\otimes1\otimes u)$. We see that $$(1\otimes u^{-1}\otimes1)(u^{-1}\otimes1\otimes1)(1\otimes u\otimes1)=\Psi_{1}(u)\otimes1=u^{-1}\otimes1\otimes1.$$ In addition, it follows from Lemma\ \ref{pan2-11} that $$u^{-1}\otimes1\otimes1=\prod_{x=1}^{n}\prod_{y=1}^{n}((a_1,a_2,a_3,x,y),(b_1,b_2,b_3,x,y)).$$
Therefore, $$\Psi_{2}(u)=(1\otimes1\otimes u^{-1})\Big(\prod_{x=1}^{n}\prod_{y=1}^{n}((a_1,a_2,a_3,x,y),(b_1,b_2,b_3,x,y))\Big)(1\otimes1\otimes u).$$
Since $1\otimes1\otimes u^{-1}=(1\otimes1\otimes u)^{-1}$, we deduce that
$$\Psi_{2}(u)=\prod_{x=1}^{n}\prod_{y=1}^{n}((a_1,a_2,u(a_3,x,y)),(b_1,b_2,u(b_3,x,y))).$$
Moreover, $\{a_1, b_1\} \cap \{a_3, b_3\}= \emptyset$ implies that $u(a_3,x,y)=(a_3,x,y)$ and $u(b_3,x,y)=(b_3,x,y)$ for any $x,y\in[n]$. Therefore,
$$\Psi_{2}(u)=\prod_{x=1}^{n}\prod_{y=1}^{n}((a_1,a_2,a_3,x,y),(b_1,b_2,b_3,x,y))=u^{-1}\otimes1\otimes1=\Psi_{1}(u)\otimes1.$$

So far, we have seen that each equation (2.1) holds for $u$, as desired.  \qed

\begin{lem}\label{pan2-3}\normalfont
Let $u=((a_1,a_2,a_3),(b_1,b_2,b_3))$ be a transposition in $S([n]^3)$ with $a_1=a_3=b_1=b_3$ and $a_1\neq a_2\neq b_2\neq b_1$. Then $u$ is stable of rank $1$.
\end{lem}
\demo One easily checks that $\Psi_{1}(u)=\Psi_{0}(u)\otimes1$. From Lemma\ \ref{pan2-2}, we see that $$\Psi_{2}(u)=\prod_{x=1}^{n}\prod_{y=1}^{n}((a_1,a_2,u(a_3,x,y)),(b_1,b_2,u(b_3,x,y))).$$
Note that $u(a_3,x,y)=(a_3,x,y)$ and $u(b_3,x,y)=(b_3,x,y)$ for any $x,y\in[n]$ except $x=a_2,y=a_3$ and $x=b_2,y=b_3$. However, we see that $u(a_3,a_2,a_3)=(a_3,b_2,b_3)$ and $u(b_3,a_2,a_3)=(b_3,b_2,b_3)$, and $u(a_3,b_2,b_3)=(a_3,a_2,a_3)$ and $u(b_3,b_2,b_3)=(b_3,a_2,a_3)$. Thus, $\Psi_{2}(u)=\Psi_{1}(u)\otimes1$ still holds. The proof of this lemma is completed.  \qed

Applying Lemma\ \ref{pan2-2} and Lemma\ \ref{pan2-3}, we see that ${\bf{(3)}}$ implies ${\bf{(2)}}$. Additionally, it is obvious that ${\bf{(2)}}$ implies ${\bf{(1)}}$. Next we will prove that ${\bf{(1)}}$ implies ${\bf{(3)}}$ by contradiction. In other words, we shall show that except for $a_1=a_3=b_1=b_3$ and $a_1\neq a_2\neq b_2\neq b_1$, if $\{a_1, b_1\} \cap \{a_3, b_3\}\neq \emptyset$ or $\{a_1, a_2\}=\{b_2, b_3\}$ or $\{b_1, b_2\}=\{a_2, a_3\}$ then $((a_1,a_2,a_3),(b_1,b_2,b_3))$ is not stable.

Consider $\{a_1, b_1\} \cap \{a_3, b_3\}\neq \emptyset$. Note that there exist four possible cases that are $a_1=a_3$ and $a_1=b_3$ and $b_1=a_3$ and $b_1=b_3$. By symmetry, it suffices to consider $a_1=a_3$ and $a_1=b_3$. Next we begin with $a_1=a_3$ that can be divided into $a_1=a_3=a_2$ and $a_1=a_3\neq a_2$.

\begin{lem}\label{pan2-4}\normalfont
Let $u=((a_1,a_2,a_3),(b_1,b_2,b_3))$ be a transposition in $S([n]^3)$ with $a_1=a_2=a_3$. Then $u$ is not stable.
\end{lem}
\demo Assume that $a_1=a_2=a_3=a$. If $a\not\in\{b_1,b_2,b_3\}$, then we note that for any $r\in N$, $$\Psi_{3r-2}(u)(\underbrace{a,a,...,a}_{3r},a)=(a,a,...,a)$$ and $$\Psi_{3r-2}(u)(\underbrace{a,a,...,a}_{3r},b_1)=(b_1,x_1...,x_{3r}).$$ Clearly, $\Psi_{3r-2}(u)\not\in S([n]^{3r})\otimes1$ as $a\neq b_1$. It follows from Lemma\ \ref{pan2-111} that $u$ is not stable. Next we divide into two subcases to discuss the case that $a\in\{b_1,b_2,b_3\}$.

{\bf{Case 1}}: $a=b_1$.

 Suppose $a\not\in\{b_2,b_3\}$. In this case, we see that for any even integer $r\in N$, $$\Psi_{3r-2}(u)(\underbrace{a,a,...,a}_{3r},a)=(b_1,\underbrace{b_2,b_3,b_2,b_3,...,b_2,b_3}_{3r}).$$ Since $a\neq b_3$, it follows that $\Psi_{3r-2}(u)\not\in S([n]^{3r})\otimes1$, and thus $u$ is not stable by Lemma\ \ref{pan2-111}.

Suppose $a=b_2$ and $a\neq b_3$. In this case, we see that for any $r\in N$, $$\Psi_{3r-2}(u)(\underbrace{a,a,...,a}_{3r},a)=(b_1,b_2,\underbrace{b_3,b_3,...,b_3}_{3r-1}).$$ Since $a\neq b_3$, we have $\Psi_{3r-2}(u)\not\in S([n]^{3r})\otimes1$, and thus $u$ is not stable.

Suppose $a=b_3$ and $a\neq b_2$. In this case, we see that for any even integer $r\in N$, $$\Psi_{3r-2}(u)(\underbrace{a,a,...,a}_{3r},a)=(b_1,b_2,\underbrace{a,a,...,a}_{3r-1})$$ and $$\Psi_{3r-2}(u)(\underbrace{a,a,...,a}_{3r},b_2)=(\underbrace{a,a,...,a}_{3r},b_2).$$Note that $\Psi_{3r-2}(u)\not\in S([n]^{3r})\otimes1$. Thus, $u$ is not stable by Lemma\ \ref{pan2-111}.

{\bf{Case 2}}: $a\neq b_1$.

In this case, we see that for any $r\in N$, $$\Psi_{3r-2}(u)(\underbrace{a,a,...,a}_{3r},a)=(a,a,...,a)$$ and $$\Psi_{3r-2}(u)(\underbrace{a,a,...,a}_{3r},b_1)=(b_1,x_1,...,x_{3r}).$$ Note $\Psi_{3r-2}(u)\not\in S([n]^{3r})\otimes1$. Thus, $u$ is not stable due to Lemma\ \ref{pan2-111}.

According to above discussions, we deduce this lemma.  \qed

\begin{lem}\label{pan2-5}\normalfont
Let $u=((a_1,a_2,a_3),(b_1,b_2,b_3))$ be a transposition in $S([n]^3)$ with $a_1=a_3\neq a_2$. Then $u$ is not stable.
\end{lem}
\demo Assume that $\{a_1,a_2,a_3\}\cap\{b_1,b_2,b_3\}=\emptyset$. We note that for any odd integer $r\in N$, $$\Psi_{2r}(u)(a_1,a_2,a_3\underbrace{a_2,a_3,...,a_2,a_3}_{2r})=(a_1,a_2,a_3\underbrace{a_2,a_3,...,a_2,a_3}_{2r})$$ and $$\Psi_{2r}(u)(a_1,a_2,a_3\underbrace{a_2,a_3,...,a_2,a_2}_{2r})=(b_1,b_2,b_3\underbrace{a_2,a_3,...,a_2,a_2}_{2r}).$$
 Clearly, $\Psi_{2r}(u)\not\in S([n]^{2r+2})\otimes1$. Then by Lemma\ \ref{pan2-111} we see that $u$ is not stable. Next we divide into two subcases to discuss the case that $\{a_1,a_2,a_3\}\cap\{b_1,b_2,b_3\}\neq\emptyset$.

{\bf{Case 1}}: $a_1=b_1$.

Suppose $b_1\not\in\{b_2,b_3\}$. In this case, for any $r\in N$, we have $$\Psi_{2r}(u)(a_1,a_2,a_3\underbrace{a_2,a_3,...,a_2,a_3}_{2r})=(b_1,b_2,b_3\underbrace{b_2,b_3,...,b_2,b_3}_{2r}).$$ Since $a_3\neq b_3$, we have $\Psi_{2r}(u)\not\in S([n]^{2r+2})\otimes1$, and thus $u$ is not stable.

Suppose $b_1\in\{b_2,b_3\}$. Since except for $a_1=a_3=b_1=b_3$ and $a_1\neq a_2\neq b_2\neq b_1$, we have $b_1=b_2\neq b_3$. In this case, we see that for any $r\in N$, $$\Psi_{2r}(u)(a_1,a_2,a_3\underbrace{a_2,a_3,...,a_2,a_3}_{2r})=(b_1,b_2,b_3\underbrace{b_2,b_3,...,b_2,b_3}_{2r}).$$ Clearly $\Psi_{2r}(u)\not\in S([n]^{2r+2})\otimes1$, and thus $u$ is not stable.

{\bf{Case 2}}: $a_1\neq b_1$.

By Lemma\ \ref{pan2-4}, we can assume that $b_1=b_2=b_3$ does not occur. Thus, we note that for any odd integer $r\in N$, $$\Psi_{2r}(u)(a_1,a_2,a_3\underbrace{a_2,a_3,...,a_2,a_3}_{2r})=(a_1,a_2,a_3\underbrace{a_2,a_3,...,a_2,a_3}_{2r})$$ and $$\Psi_{2r}(u)(a_1,a_2,a_3\underbrace{a_2,a_3,...,a_2,a_2}_{2r})=(b_1,b_2,b_3\underbrace{a_2,a_3,...,a_2,a_2}_{2r}).$$
Clearly, $\Psi_{2r}(u)\not\in S([n]^{2r+2})\otimes1$. Then by Lemma\ \ref{pan2-111} we see that $u$ is not stable.

According to above discussions, we deduce this lemma.  \qed

So far, we have proved the case that $a_1=a_3$. By symmetry, it suffices to consider the case that $a_1=b_3$ and $a_1\neq a_3$ and $b_1\neq b_3$.

\begin{lem}\label{pan2-6}\normalfont
Let $u=((a_1,a_2,a_3),(b_1,b_2,b_3))$ be a transposition in $S([n]^3)$ with $a_1=b_3$ and $a_1\neq a_3$ and $b_1\neq b_3$. Then $u$ is not stable.
\end{lem}
\demo Since $a_1=b_3$ and $a_1\neq a_3$ and $b_1\neq b_3$, it follows that for any $r\in N$, $$\Psi_{2r}(u)(\underbrace{b_1,b_2,...,b_1,b_2}_{2r},b_1,b_2,b_3)=(a_1,\underbrace{a_2,a_3,...,a_2,a_3}_{2r}).$$ Since $a_3\neq b_3$, it follows that $\Psi_{2r}(u)\not\in S([n]^{2r+2})\otimes1$. Thus, $u$ is not stable by Lemma\ \ref{pan2-111}, as desired.  \qed

Up to now, we have proved that if $\{a_1, b_1\} \cap \{a_3, b_3\}\neq \emptyset$ then $u$ is not stable. So we assume that $\{a_1, b_1\} \cap \{a_3, b_3\}=\emptyset$. In this case, if $\{a_1, a_2\}=\{b_2, b_3\}$ then $b_2=a_1$ and $b_3=a_2$, and if $\{b_1, b_2\}=\{a_2, a_3\}$ then $a_2=b_1$ and $a_3=b_2$. Next we start to prove the case that $a_2=b_1$ and $a_3=b_2$. Then by symmetry, the other case can be deduced immediately.

\begin{lem}\label{pan2-7}\normalfont
Let $u=((a_1,a_2,a_3),(b_1,b_2,b_3))$ be a transposition in $S([n]^3)$ with $a_2=b_1$ and $a_3=b_2$ and $\{a_1, b_1\} \cap \{a_3, b_3\}=\emptyset$. Then $u$ is not stable.
\end{lem}
\demo Note $b_1\neq b_3$. For any $r\in N$, we see that $$\Psi_{3r-2}(u)(\underbrace{a_1,a_1,...,a_1}_{3r-3},a_1,a_2,a_3,b_1)=(b_1,b_2,\underbrace{b_3,b_3,...,b_3}_{3r-2},b_1)$$ and $$\Psi_{3r-2}(u)(\underbrace{a_1,a_1,...,a_1}_{3r-3},a_1,a_2,a_3,b_3)=(\underbrace{a_1,a_1,...,a_1}_{3r-2},b_1,b_2,b_3).$$ Thus, $\Psi_{3r-2}(u)\not\in S([n]^{3r})\otimes1$. By Lemma\ \ref{pan2-111}, we see that $u$ is not stable, as desired.  \qed

\begin{lem}\label{pan2-8}\normalfont
Let $u=((a_1,a_2,a_3),(b_1,b_2,b_3))$ be a transposition in $S([n]^3)$ with $b_2=a_1$ and $b_3=a_2$. Then $u$ is not stable.
\end{lem}

Up to now we have completed the proof of Theorem\ \ref{pan1-2}.



\end{document}